\documentclass[11pt]{article}
\usepackage[T1]{fontenc}
\usepackage{amsmath,amssymb,amsthm}
\usepackage{url}
\usepackage{hyperref}
\usepackage[margin=1.15in]{geometry}
\usepackage{xcolor}
\theoremstyle{plain}
\newtheorem{theorem}{Theorem}[section]

\newtheorem{proposition}[theorem]{Proposition}

\newtheorem{conjecture}[theorem]{Conjecture}
\theoremstyle{definition}
\newtheorem{definition}[theorem]{Definition}
\newtheorem{problem}[theorem]{Problem}
\newtheorem{question}[theorem]{Question}
\newtheorem{nonitalicexample}[theorem]{Example}
\theoremstyle{remark}
\newtheorem{remark}[theorem]{Remark}
\newcommand{\PP}{\mathbb P}
\newcommand{\Aaff}{\mathbb A}
\newcommand{\CC}{\mathbb C}
\newcommand{\KK}{\mathbb K}
\newcommand{\OO}{\mathcal O}
\newcommand{\cJ}{\mathcal J}
\newcommand{\fraka}{\mathfrak a}
\newcommand{\frakb}{\mathfrak b} 

\newcommand{\bR}{\mathbf R}   
\DeclareMathOperator{\ord}{ord}    
\begin{document}

\title{Containment problems for fat point schemes via asymptotic multiplier ideals}
\author{V\'ictor Gonz\'alez-Alonso, Piotr Pokora}
\date{\today}
\maketitle
\thispagestyle{empty}
\begin{abstract}
We study a multiplier-ideal approach to the containment problem for radical ideals of points in projective space. The expected asymptotic multiplier-ideal computation depends on controlling residual base loci of the relevant symbolic linear systems. We isolate this issue in the notion of clean asymptotic multiplier ideals. Under this hypothesis, a failure of
\[
        I^{(2r-1)}\not\subseteq I^r
\]
forces
\[
        2r\widehat{\alpha}(I)\geq \alpha(I^{(2r-1)})+3,
\]
and in particular \(\widehat{\alpha}(I)\geq 3\). We show that the dual-Hesse configuration is clean although its symbolic systems have a genuinely non-trivial residual asymptotic base part. We also verify cleanliness for star configurations of lines, giving a further infinite family of non-monomial examples, and formulate a cleanliness conjecture for singular loci of line arrangements.
\end{abstract}

\noindent
\textbf{Keywords:} containment problems, line configurations, multiplier ideals, symbolic powers of homogeneous ideals, Waldschmidt constants.

\medskip

\noindent
\textbf{Mathematics Subject Classification (2020):} 14C20, 14F18, 14J26, 14N20, 13A15, 13F20.

\section{Preliminaries}

In this paper we consider the planar case of containment problems. For a general introduction to the subject we refer to the survey \cite{SzSz}. Let
\[
        \mathcal P=\{P_1,\ldots,P_s\}
\]
be a finite set of mutually distinct points in the projective plane over an algebraically closed field \(\KK\) of characteristic \(0\). We denote by
\[
        I:=I(\mathcal P)=I(P_1)\cap \cdots \cap I(P_s)\subset \KK[x,y,z]
\]
the corresponding homogeneous radical ideal.

For \(m\geq 1\), the \(m\)-th symbolic power of \(I\) is given by
\[
        I^{(m)}
        =
        I(P_1)^m\cap \cdots \cap I(P_s)^m.
\]
Equivalently, \(I^{(m)}\) is the ideal generated by homogeneous forms vanishing to order at least \(m\) at every point \(P_i\). This description agrees with the original definition of symbolic powers by the Nagata--Zariski theorem.

The inclusion
\[
        I^m\subseteq I^{(m)}
\]
is immediate. In general, however, this inclusion is far from being an equality, and one is led to ask for reverse containments with different exponents.

\begin{problem}
Decide for which pairs of positive integers \((m,r)\) the following containment holds:
\[
        I^{(m)}\subseteq I^r.
\]
\end{problem}

A fundamental uniform result was proved by Ein, Lazarsfeld, and Smith in characteristic \(0\) \cite{ELS}, and by Hochster and Huneke in positive characteristic \cite{HH}.

\begin{theorem}\label{thm:ELS}
Let \(\mathcal P\) be a finite set of points in \(\PP^n\), and let \(I=I(\mathcal P)\) be the associated radical ideal. Then, for every \(r\geq 1\), one has
\[
        I^{(nr)}\subseteq I^r.
\]
\end{theorem}

Huneke asked whether this uniform containment is sharp. In the planar case this leads to the following problem, often attributed to Bocci, Harbourne, and Huneke; see \cite[Problem 1.5]{SzSz}.

\begin{problem}\label{prob:BHH}
Let \(\mathcal P\subset \PP^2\) be a finite set of points, and let \(I=I(\mathcal P)\). Is it true that, for every \(r\geq 1\), one has
\begin{equation}\label{eq:BHH}
        I^{(2r-1)}\subseteq I^r?
\end{equation}
\end{problem}

The first particularly important case is
\begin{equation}\label{eq:I3I2}
        I^{(3)}\subseteq I^2.
\end{equation}
We recall that Dumnicki, Szemberg, and Tutaj-Gasi\'nska \cite{DST} constructed the first counterexample to \eqref{eq:I3I2} over the complex numbers. Their example is given by the singular locus of the dual-Hesse configuration of \(9\) lines and \(12\) triple points. Most known counterexamples to \eqref{eq:I3I2} also arise from line configurations. For example, Klein's and Wiman's configurations give further counterexamples; see \cite{BNC,Wal,Sec}.

It remains difficult to find counterexamples to \eqref{eq:I3I2}. One useful criterion is due to Bocci and Harbourne. They showed, in particular, that if \(I\subset \KK[x,y,z]\) is a homogeneous radical ideal, then
\[
        I^{(2r-1)}_d\subseteq I^{(r)}_d=I^r_d
\]
for every \(d\geq \operatorname{reg}(I^r)\); see \cite[Lemma 2.3.3]{BH}. Thus, in order to find a counterexample, one first looks for a form in \(I^{(3)}\) of degree smaller than \(\operatorname{reg}(I^2)\). In the dual-Hesse example, the product of the \(9\) lines has degree \(9\), while \(\operatorname{reg}(I^2)=10\).

The aim of this note is to revisit a possible multiplier-ideal approach to numerical obstructions for the failure of \eqref{eq:BHH}. We will explain that the most direct computation of the relevant asymptotic multiplier ideals is not automatic: it requires control of residual base loci on the blow-up of the point set. Under this additional hypothesis, one obtains a useful lower bound for the Waldschmidt constant.

We recall the basic definitions that we will use in this note.

\begin{definition}
Let \(I\subset \KK[x_0,\ldots,x_n]\) be a homogeneous ideal. The initial degree of \(I\) is
\[
        \alpha(I)=\min\{t\mid I_t\neq 0\}.
\]
\end{definition}

\begin{definition}
Let \(I\subset \KK[x_0,\ldots,x_n]\) be a homogeneous radical ideal. The Waldschmidt constant of \(I\) is
\[
        \widehat{\alpha}(I)
        =
        \inf_{m\geq 1}\frac{\alpha(I^{(m)})}{m}
        =
        \lim_{m\to\infty}\frac{\alpha(I^{(m)})}{m}.
\]
\end{definition}
\noindent
For brevity, we will write
\[
        \alpha_m=\alpha(I^{(m)}),
        \qquad
        \widehat{\alpha}=\widehat{\alpha}(I).
\]

\section{Clean asymptotic multiplier ideals}

For the multiplier-ideal part we write the ground field as \(\CC\), in accordance with the standard references. All arguments in this section are algebraic and remain valid over an algebraically closed field of characteristic \(0\).

We use standard notation and conventions for multiplier ideals; see, for instance, \cite[Chapter~9]{LazarsfeldPositivityII}. Let \(X\) be a smooth variety and let
\[
        \fraka\subset \OO_X
\]
be a coherent ideal sheaf. A log resolution of \(\fraka\) is a proper birational morphism
\[
        \pi:X'\longrightarrow X
\]
such that
\[
        \fraka\OO_{X'}=\OO_{X'}(-F)
\]
for an effective divisor \(F\subset X'\), and such that \(F+\operatorname{Exc}(\pi)\) has simple normal crossings.

\begin{definition}
With the notation above, the multiplier ideal of \(\fraka\) with exponent \(c\in \mathbb R_{\geq 0}\) is
\[
        \cJ(\fraka^c)
        =
        \cJ(c\cdot \fraka)
        =
        \pi_*\OO_{X'}(K_{X'/X}-\lfloor cF\rfloor).
\]
\end{definition}
\begin{remark}
The definition is independent of the chosen log resolution.
\end{remark}
Let \(\{\fraka_d\}_{d\geq 1}\) be a multiplicative graded system of coherent ideal sheaves, so that
\[
        \fraka_d\fraka_e\subseteq \fraka_{d+e}
\]
for all \(d,e\geq 1\). For fixed \(c\geq 0\) and \(d\geq 1\), the family
\[
        \left\{\cJ\left(\frac{c}{N}\cdot \fraka_{Nd}\right)\right\}_{N\geq 1}
\]
has a unique maximal member. Equivalently, for all sufficiently divisible \(N\), the corresponding multiplier ideals agree with this maximal member. We denote it by
\[
        \cJ(c\cdot \|\fraka_d\|).
\]
This is the asymptotic multiplier ideal associated to the graded system; see \cite{ELS,Teitler}. We will mostly use \(c=1\) and write simply \(\cJ(\|\fraka_d\|)\). Thus no literal limit or union of ideals is intended by the double-bar notation.

The following properties are standard; see \cite[Proposition 1.7]{ELS}:
\[
        \fraka_d\subseteq \cJ(\|\fraka_d\|)
\]
and
\[
        \cJ(\|\fraka_{md}\|)
        \subseteq
        \cJ(\|\fraka_d\|)^m.
\]
In particular, if
\[
        \cJ(\|\fraka_d\|)\subseteq \fraka_1,
\]
then
\[
        \fraka_{md}\subseteq \fraka_1^m
\]
for every \(m\geq 1\). This is the mechanism behind the proof of Theorem \ref{thm:ELS}.

We now specialize to symbolic powers of point ideals. Let
\[
        I=I(P_1,\ldots,P_s)\subset \CC[x_0,\ldots,x_n]
\]
be the radical ideal of a finite set of points in \(\PP^n\). The family \(\{I^{(m)}\}_{m\geq 1}\) is a multiplicative graded system.

Passing to affine cones, \(I\) defines the union of lines
\[
        \ell_1,\ldots,\ell_s\subset \Aaff^{n+1}
\]
through the origin. Let
\begin{equation}\label{eq:pi}
        \pi:X'\longrightarrow \Aaff^{n+1}
\end{equation}
be the morphism obtained by first blowing up the origin and then blowing up the strict transforms of the lines \(\ell_i\). Let \(F_0\) be the strict transform of the exceptional divisor over the origin, and for each \(i=1,\ldots,s\) let \(F_i\) be the exceptional divisor over \(\ell_i\).

The divisors \(F_0,F_1,\ldots,F_s\) record the basic valuations of the symbolic powers. Namely, every form in \(I^{(m)}\) vanishes at the origin to order at least \(\alpha_m\), and along every line \(\ell_i\) to order at least \(m\). In other words,
\[
        I^{(m)}\OO_{X'}\subseteq \OO_{X'}\bigl(-\alpha_mF_0-m(F_1+\cdots+F_s)\bigr),
\]
so that
\begin{equation}\label{eq:residual-ideal}
        \frakb_m:=I^{(m)}\OO_{X'}\cdot \OO_{X'}\bigl(\alpha_mF_0+m(F_1+\cdots+F_s)\bigr)\subseteq \OO_{X'}
\end{equation}
is an ideal sheaf on \(X'\), which we call the \emph{residual ideal} of \(I^{(m)}\). If
\[
        \mu_m:Y_m\longrightarrow \Aaff^{n+1}
\]
is a log resolution of \(I^{(m)}\) dominating \(X'\), say \(\mu_m=\pi\circ\sigma_m\), and
\[
        I^{(m)}\OO_{Y_m}=\OO_{Y_m}(-D_m),
\]
then
\begin{equation}\label{eq:Dm-lower-bound}
        D_m
        \geq
        \alpha_mE_0+m(E_1+\cdots+E_s),
\end{equation}
where \(E_i:=\sigma_m^*F_i\) denotes the total transform of \(F_i\) on \(Y_m\), for \(i=0,\ldots,s\). We therefore write
\begin{equation}\label{eq:Dm-with-residual}
        D_m
        =
        \alpha_mE_0
        +
        m(E_1+\cdots+E_s)
        +
        R_m,
\end{equation}
where \(R_m\) is an effective residual contribution; namely, \(\frakb_m\OO_{Y_m}=\OO_{Y_m}(-R_m)\). It records divisorial or higher-codimensional base-locus phenomena in the linear system defined by \(I^{(m)}\) not detected by the expected valuations along the divisors \(F_0,F_1,\ldots,F_s\).

\begin{remark}\label{rem:incorrect-equality}
In general \(\pi\) as in \eqref{eq:pi} is not a log-resolution of \(I^{(m)}\), i.e.\ there is no equality
\[
        I^{(m)}\OO_{X'}
        =
        \OO_{X'}\bigl(-\alpha_mF_0-m(F_1+\cdots+F_s)\bigr),
\]
or equivalently \(\frakb_m\neq\OO_{X'}\). The dual-Hesse example below (Proposition~\ref{prop:dual-hesse-clean}) shows explicitly that this can fail even for a very structured line arrangement: to obtain a log-resolution it is still necessary to blow up the strict transforms of the lines defining the arrangement, and \(R_m\) is non-zero along the additional exceptional divisors (at least when \(3\mid m\)). 
Nevertheless, these extra conditions imposed by \(R_m\) turn out to be harmless for the asymptotic multiplier ideal.
\end{remark}

\begin{remark}
It is also not clear whether all symbolic powers \(I^{(m)}\) have a common log-resolution at all. This is the case for certain line arrangements (see the proofs of Propositions~\ref{prop:dual-hesse-clean} and \ref{prop:star-clean}), but this might fail for more general point configurations.
\end{remark}

Since a common log-resolution for all symbolic powers \(I^{(m)}\) might not exist, it is convenient to formulate the additional conditions imposed by the \(R_m\) in terms of valuations. If \(v\) is a valuation of \(\CC(x_0,\ldots,x_n)\) with centre on \(\Aaff^{n+1}\) (for instance, a divisorial valuation over \(\Aaff^{n+1}\)), set
\begin{equation}\label{eq:asymptotic-order}
        v(\|I^{(m)}\|)
        :=
        \lim_{N\to\infty}\frac{v(I^{(Nm)})}{N}
        =
        \inf_{N\geq 1}\frac{v(I^{(Nm)})}{N}.
\end{equation}
The limit exists by subadditivity (Fekete's lemma), since \(I^{(k)}I^{(l)}\subseteq I^{(k+l)}\) gives \(v(I^{(k+l)})\leq v(I^{(k)})+v(I^{(l)})\). For the same reason
\[
        v(\|I^{(m)}\|)=m\,v(\|I^{(\bullet)}\|),
        \qquad\text{where}\quad
        v(\|I^{(\bullet)}\|):=\lim_{k\to\infty}\frac{v(I^{(k)})}{k}=\inf_{k\geq 1}\frac{v(I^{(k)})}{k},
\]
so all asymptotic orders are determined by the normalised order \(v(\|I^{(\bullet)}\|)\). In particular,
\[
        v_{F_0}(\|I^{(m)}\|)=m\widehat{\alpha}(I),
        \qquad
        v_{F_i}(\|I^{(m)}\|)=m.
\]

For such a valuation \(v\), let \(v(F_i)\) denote the value of \(v\) on a local equation of \(F_i\) at the centre of \(v\) on \(X'\); thus \(v(F_i)=0\) if this centre is not contained in \(F_i\). By \eqref{eq:residual-ideal},
\[
        v(I^{(k)})=\alpha_k\,v(F_0)+k\sum_{i=1}^s v(F_i)+v(\frakb_k)
        \qquad\text{for every }k\geq 1.
\]
Dividing by \(k\) and letting \(k\to\infty\) gives
\begin{equation}\label{eq:residual-order}
        v(\|I^{(\bullet)}\|)
        =
        \widehat{\alpha}(I)\,v(F_0)+\sum_{i=1}^s v(F_i)+r(v),
        \qquad
        r(v):=\lim_{k\to\infty}\frac{v(\frakb_k)}{k}\geq 0 .
\end{equation}

\begin{definition}\label{def:residual-b-divisor}
We call \(r(v)\) the \emph{residual asymptotic order} of \(v\). The \emph{residual asymptotic base b-divisor} of the symbolic graded system \(\{I^{(m)}\}_{m\geq 1}\) is the \(\mathbb R\)-b-divisor \(\bR\) over \(X'\) whose coefficient along a prime divisor \(E\) over \(X'\) is \(r(\ord_E)\).
\end{definition}

In other words, if each \(R_k\) is regarded as the b-divisor of the residual ideal \(\frakb_k\) (its trace on \(Y_k\) is the divisor \(R_k\) of \eqref{eq:Dm-with-residual}), then
\[
        \bR=\lim_{k\to\infty}\frac1k R_k
\]
coefficientwise. Note that this is the limit of the normalised divisors \(R_k/k\), not of the \(R_k\) themselves. By the homogeneity of asymptotic orders, the residual asymptotic base b-divisor of the graded system \(\{I^{(Nm)}\}_{N\geq 1}\) is simply \(m\bR\), so a single object governs all \(m\) at once. The individual divisors \(R_m\) carry additional, non-asymptotic information, which is invisible to the asymptotic multiplier ideals. Note moreover that the trace of $\bR$ on $X'$ itself is zero.

If \(r(v)>0\) for some divisorial valuation \(v\), then \(v(\frakb_k)>0\) for all \(k\gg 0\). Hence \(\pi\) is not a log-resolution of \(I^{(k)}\) for \(k\gg0\), and \(v\) contributes to \(R_k\) on every log-resolution \(Y_k\) on which it is realised by a divisor. It is \(r(v)\), and not \(v(\|I^{(m)}\|)\), that detects this: \(v(\|I^{(m)}\|)>0\) for every valuation whose centre on \(\Aaff^{n+1}\) lies in \(\ell_1\cup\cdots\cup\ell_s\).

With these notations, the asymptotic multiplier ideals can be described valuatively. For a prime divisor \(E\) over \(\Aaff^{n+1}\), i.e.\ a prime divisor on some normal birational model \(\mu\colon Y\to\Aaff^{n+1}\), put
\[
        k_E:=\ord_E\bigl(K_{Y/\Aaff^{n+1}}\bigr),
        \qquad
        A(E):=k_E+1 .
\]
Neither depends on the choice of \(Y\). We also write \(k_v=k_E\) and \(A(v)=A(E)\) for \(v=\ord_E\). In standard terminology, \(k_E\) is the \emph{discrepancy} and \(A(E)\) the \emph{log discrepancy} of \(E\) over \(\Aaff^{n+1}\); see \cite[Section~2.3]{KM}. The notation \(A(\cdot)\) is the one used in \cite{JM}, where the log discrepancy is extended to arbitrary real valuations.

Compute \(\cJ(\frac1N\cdot I^{(Nm)})\) on a log resolution on which \(E\) appears, and use
\[
        \frac{\ord_E(I^{(Nm)})}{N}\geq \ord_E(\|I^{(m)}\|).
\]
This gives
\begin{equation}\label{eq:valuative}
        \cJ(\|I^{(m)}\|)
        \subseteq
        \Bigl\{f\in\CC[x_0,\ldots,x_n]\ \Big|\
        \ord_E(f)\geq \bigl\lfloor \ord_E(\|I^{(m)}\|)\bigr\rfloor-k_E
        \ \text{for all } E\Bigr\},
\end{equation}
where \(E\) runs through all prime divisors over \(\Aaff^{n+1}\). Conversely, by \cite{JM}, \(f\in\cJ(\|I^{(m)}\|)\) if and only if
\begin{equation}\label{eq:valuative-JM}
        v(f)+A(v)>v(\|I^{(m)}\|)
\end{equation}
for every real valuation \(v\) with centre on \(\Aaff^{n+1}\) and \(A(v)<\infty\). For a divisorial valuation \(v=\ord_E\), \eqref{eq:valuative-JM} is equivalent to the condition in \eqref{eq:valuative}, because \(\ord_E(f)\) and \(k_E\) are integers. In general, however, the strict inequalities \eqref{eq:valuative-JM} must be tested on all (quasi-monomial) valuations, not only on divisorial ones.

This motivates the following definition.

\begin{definition}\label{def:clean}
Let \(I=I(P_1,\ldots,P_s)\subset \CC[x_0,\ldots,x_n]\). We say that \(I\) has \emph{clean asymptotic multiplier ideals} if, for every \(m\geq n\),
\begin{equation}\label{eq:clean-formula}
        \cJ(\|I^{(m)}\|)
        =
        \pi_*
        \OO_{X'}
        \left(
        -(\lfloor m\widehat{\alpha}\rfloor-n)F_0
        -(m-n+1)(F_1+\cdots+F_s)
        \right).
\end{equation}
\end{definition}
Since \(k_{F_0}=n\) and \(k_{F_i}=n-1\), the right-hand side of \eqref{eq:clean-formula} is exactly the ideal cut out by the conditions in \eqref{eq:valuative} for \(E=F_0,F_1,\ldots,F_s\) alone (cf.\ \eqref{eq:pushforward} below). In particular, the inclusion \(\subseteq\) in \eqref{eq:clean-formula} always holds. Thus \(I\) is clean if and only if the valuations other than \(\ord_{F_0},\ldots,\ord_{F_s}\) impose no additional condition, i.e.\ if and only if every \(f\) in the right-hand side of \eqref{eq:clean-formula} satisfies \eqref{eq:valuative-JM} for every \(v\). By \eqref{eq:residual-order}, the inequality \eqref{eq:valuative-JM} reads
\[
        v(f)+A(v)
        >
        m\Bigl(\widehat{\alpha}\,v(F_0)+\sum_{i=1}^s v(F_i)\Bigr)+m\,r(v).
\]
So cleanliness means that the residual term \(m\,r(v)\) coming from \(\bR\) is never decisive for elements that satisfy the expected conditions along \(F_0,\ldots,F_s\).

We stress that this is a \emph{relative} redundancy condition. It is much weaker than requiring each additional valuation to impose no condition at all, i.e.\ \(v(\|I^{(m)}\|)<A(v)\) (for divisorial \(v\), equivalently \(\lfloor v(\|I^{(m)}\|)\rfloor\leq k_v\)). For instance, the divisors \(G_j\) in the proof of Proposition~\ref{prop:dual-hesse-clean} satisfy \(v_{G_j}(\|I^{(m)}\|)=10m/3>4=A(G_j)\) for all \(m\geq 2\). Nevertheless, the conditions they impose are implied by those along \(F_0,F_1,\ldots,F_{12}\).

We present now an explicit example where the residual divisors \(R_m\) are genuinely present but impose only redundant conditions on the asymptotic multiplier ideals.

\begin{proposition}\label{prop:dual-hesse-clean}
Let
\[
        Q=(x^3-y^3)(y^3-z^3)(z^3-x^3)
\]
be the defining equation of the dual-Hesse arrangement, let \(Z=\operatorname{Sing}(Q)\) be its twelve triple points, and put
\[
        I=I(Z)\subset \CC[x,y,z].
\]
Then the residual asymptotic base b-divisor \(\bR\) of the symbolic graded system is non-zero. Its coefficient equals \(\frac13\) along the exceptional divisor obtained by blowing up the strict transform, in \(F_0\cong\operatorname{Bl}_Z\PP^2\), of each of the nine Hesse lines. Nevertheless, \(I\) has clean asymptotic multiplier ideals. More precisely, for every \(m\geq 2\),
\begin{equation}\label{eq:hesse-J}
        \cJ(\|I^{(m)}\|)
        =
        \pi_*\OO_{X'}\bigl(-(3m-2)F_0-(m-1)(F_1+\cdots+F_{12})\bigr)
        =
        \bigoplus_{d\geq 3m-2}I^{(m-1)}_d .
\end{equation}
\end{proposition}

\begin{proof}
Let
\[
        \rho : X=\operatorname{Bl}_Z\PP^2\longrightarrow \PP^2
\]
be the blow-up of the twelve triple points. Denote by \(H\) the pull-back of a line and by \(E_1,\ldots,E_{12}\) the exceptional divisors. Let \(L_1,\ldots,L_9\) be the nine lines of the dual-Hesse arrangement and for each \(j=1,\ldots,9\) let \(C_j\subset X\) denote the strict transform of \(L_j\). Each \(L_j\) contains four points of \(Z\), while each point of \(Z\) lies on three Hesse lines. Therefore
\[
        C_j\sim H-\sum_{P_i\in L_j}E_i,
        \qquad
        C_j^2=-3,
\]
the curves \(C_1,\ldots,C_9\) are pairwise disjoint on \(X\), and
\[
        \sum_{j=1}^9 C_j = 9H-3\sum_{i=1}^{12}E_i.
\]
Set
\[
        B=4H-\sum_{i=1}^{12}E_i.
\]
The divisor \(B\) is base-point-free, hence nef. Indeed, the quartics
\[
        q_1=x(y^3-z^3),\qquad q_2=y(z^3-x^3),\qquad q_3=z(x^3-y^3)
\]
have precisely \(Z\) as their common zero locus, and their tangent cones at each \(P\in Z\) consist of different lines, proving that \(B\) is base-point-free.

For \(\lambda\geq 3\), put
\[
        D_\lambda=\lambda H-\sum_{i=1}^{12}E_i.
\]
For \(3\leq \lambda\leq 4\), one has
\begin{equation}\label{eq:hesse-zariski}
        D_\lambda
        =
        (\lambda-3)B+
        \frac{4-\lambda}{3}\sum_{j=1}^9 C_j.
\end{equation}
Since \(B\cdot C_j=0\), the \(C_j\)'s are disjoint, and \(C_j^2=-3\), \eqref{eq:hesse-zariski} is the Zariski decomposition of \(D_\lambda\), with positive resp.\ negative parts given by
\[
        P_\lambda=(\lambda-3)B,
        \qquad
        N_\lambda=\frac{4-\lambda}{3}\sum_{j=1}^9 C_j.
\]
For \(\lambda\geq 4\), the equality
\[
        D_\lambda=B+(\lambda-4)H
\]
shows that \(D_\lambda\) is base-point-free. Thus the negative part of \(D_\lambda\) is non-zero only for \(3\leq\lambda<4\): this is the promised explicit residual contribution.

From these computations we also obtain
\[
        \widehat\alpha(I)=3.
\]
Indeed, \(Q\) has degree \(9\) and multiplicity \(3\) at every point of \(Z\), so \(\widehat\alpha(I)\leq 3\). Conversely, if \(dH-k\sum E_i\) is effective, then, since \(B\) is nef,
\[
        0\leq \Bigl(dH-k\sum E_i\Bigr)\cdot B = 4d-12k,
\]
which gives \(d/k\geq 3\).

Now pass to the affine cone and the composition of blow-ups \(\pi:X'\to\Aaff^3\) from \eqref{eq:pi}. As above let \(F_0\) be the exceptional divisor over the origin and \(F_i\) the exceptional divisor over the cone-line corresponding to \(P_i\). By construction, there is a natural identification \(F_0= X=\operatorname{Bl}_Z\PP^2\). Let \(\mu\colon Y\to\Aaff^3\) denote the further blow-up of the nine curves \(C_j\subset F_0\), and denote the corresponding exceptional divisors by \(G_j\). By \eqref{eq:asymptotic-order},
\[
        v_{F_0}(\|I^{(m)}\|)=3m,
        \qquad
        v_{F_i}(\|I^{(m)}\|)=m.
\]
We compute the extra orders along \(G_j\). Let \(N\) be a positive integer, set \(k=Nm\), and let \(f\in I^{(k)}_d\) be a homogeneous form, with \(\lambda=d/k\geq 3\). The form \(f\) defines a curve in \(X\) with class \(dH-k\sum E_i=kD_\lambda\). If \(3\leq\lambda\leq 4\), this curve contains the negative part \(kN_\lambda\), so \(\ord_{L_j}(f)\geq k(4-\lambda)/3\). Since blowing up \(C_j\subset F_0\) adds the order along \(F_0\) (i.e.\ the degree \(d\)) to the order along \(C_j\), we have
\[
        v_{G_j}(f)\geq d+k\,\frac{4-\lambda}{3}=k\,\frac{2\lambda+4}{3}\geq \frac{10k}{3},
\]
the middle expression being increasing in \(\lambda\). If \(\lambda\geq 4\), there is no residual \(C_j\)-part and simply \(v_{G_j}(f)\geq d\geq 4k\). Since \(G_j\) is invariant under the natural \(\CC^*\)-action, \(v_{G_j}(I^{(k)})\) is attained on homogeneous forms. Hence \(v_{G_j}(I^{(Nm)})/N\geq 10m/3\) for every \(N\), and therefore
\[
        v_{G_j}(\|I^{(m)}\|)=\inf_{N\geq1}\frac{v_{G_j}(I^{(Nm)})}{N}\geq\frac{10m}{3}.
\]
This is actually an equality. If \(N\) is divisible by \(3\), then
\[
        Q^{Nm/3}\in I^{(Nm)}
\]
has degree \(3Nm\) and contains each Hesse line with multiplicity \(Nm/3\). Therefore
\[
        v_{G_j}\left(Q^{Nm/3}\right)=3Nm+\frac{Nm}{3}=\frac{10Nm}{3},
\]
so \(v_{G_j}(I^{(Nm)})\leq \frac{10Nm}{3}\), and division by \(N\) gives
\begin{equation}\label{eq:hesse-vG}
        v_{G_j}(\|I^{(m)}\|)=\frac{10m}{3}.
\end{equation}
Since \(C_j\subset F_0\) and \(C_j\not\subset F_i\), we have \(v_{G_j}(F_0)=1\) and \(v_{G_j}(F_i)=0\). By \eqref{eq:residual-order}, the residual asymptotic order is therefore
\[
        r(\ord_{G_j})=\frac{10}{3}-3=\frac13>0,
\]
which is the coefficient of \(C_j\) in the negative part \(N_3\) of \(D_3=3H-\sum E_i\). This proves \(\bR\neq 0\). In particular, if \(3\mid m\), then \(\alpha_m=3m\) and \(v_{G_j}(\frakb_m)\geq \frac{10m}{3}-3m>0\), so \(R_m\neq 0\).

With the notation of \eqref{eq:valuative}, the discrepancies of the divisors above are
\[
        k_{F_0}=2,\qquad k_{F_i}=1,\qquad k_{G_j}=3.
\]
Consequently the asymptotic multiplier ideal is cut out by
\[
        v_{F_0}(f)\geq 3m-2,
        \qquad
        v_{F_i}(f)\geq m-1,
\]
and
\begin{equation}\label{eq:hesse-extra}
        v_{G_j}(f)\geq \left\lfloor\frac{10m}{3}\right\rfloor-3 .
\end{equation}
The first two conditions define exactly
\[
        \bigoplus_{d\geq 3m-2}I^{(m-1)}_d .
\]

It remains to prove that \eqref{eq:hesse-extra} is automatic. Let \(f\in I^{(m-1)}_d\) be homogeneous, with \(d\geq 3m-2\), fix a Hesse line \(L_j\), and put
\[
        a=\ord_{L_j}(f).
\]
If \(a\geq m-1\), then
\[
        v_{G_j}(f)=d+a\geq (3m-2)+(m-1)=4m-3,
\]
which is stronger than \eqref{eq:hesse-extra}.

If \(0\leq a<m-1\), then the residual curve \(f/L_j^a\) has degree \(d-a\) and multiplicity at least \(m-1-a\) at each of the four points in \(Z\cap L_j\). B\'ezout gives immediately
\[
        d-a\geq 4(m-1-a),
\]
so
\begin{equation}\label{eq:hesse-bezout}
        a\geq\left\lceil\frac{4(m-1)-d}{3}\right\rceil .
\end{equation}
If \(d\geq 4(m-1)\), then
\[
        v_{G_j}(f)=d+a\geq d\geq 4m-4\geq \left\lfloor\frac{10m}{3}\right\rfloor-3.
\]
If \(d<4(m-1)\), then \eqref{eq:hesse-bezout} gives
\[
        v_{G_j}(f)\geq d+\left\lceil\frac{4(m-1)-d}{3}\right\rceil .
\]
The right-hand side is non-decreasing as \(d\) increases through the integers, hence its minimum in the relevant range is attained at \(d=3m-2\). Therefore
\[
        v_{G_j}(f)\geq 3m-2+\left\lceil\frac{m-2}{3}\right\rceil .
\]
For completeness, write \(m=3q+r\), \(r=0,1,2\). The last lower bound equals respectively
\[
        10q-2,\qquad 10q+1,\qquad 10q+4,
\]
whereas \(\lfloor 10m/3\rfloor-3\) equals respectively
\[
        10q-3,\qquad 10q,\qquad 10q+3.
\]
Thus \eqref{eq:hesse-extra} always holds.

Finally, there are no further asymptotic residual valuations. Indeed, after removing the negative part in \eqref{eq:hesse-zariski}, the positive part is either \((\lambda-3)B\) or \(B+(\lambda-4)H\), hence base-point-free. Therefore no further base points, including infinitely near ones, remain asymptotically. The \(G_j\)-conditions are redundant, so \eqref{eq:hesse-J} follows. This proves that \(I\) is clean despite its non-zero residual asymptotic base b-divisor \(\bR\).
\end{proof}

We next record the elementary monomial benchmark.

\begin{nonitalicexample}\label{ex:coordinate-triangle-clean}
Let
\[
        \mathcal A=V(xyz)=\{x=0\}\cup\{y=0\}\cup\{z=0\}\subset \PP^2
\]
be the fundamental triangle, and let
\[
        P=\operatorname{Sing}(\mathcal A)=\{[1:0:0],[0:1:0],[0:0:1]\}.
\]
Then
\[
        I=I(P)=(xy,xz,yz)\subset \CC[x,y,z].
\]
For every \(m\geq 1\),
\[
        I^{(m)}=(x,y)^m\cap (x,z)^m\cap (y,z)^m,
\]
so a monomial \(x^ay^bz^c\) belongs to \(I^{(m)}\) precisely when
\[
        a+b\geq m,\qquad a+c\geq m,\qquad b+c\geq m.
\]
By Howald's formula \cite{HowaldMonomialMultiplierIdeals}, for sufficiently divisible \(N\),
\[
        x^ay^bz^c\in \cJ\left(\frac1N\cdot I^{(Nm)}\right)
\]
if and only if
\[
        (a+1)+(b+1)>m,\qquad (a+1)+(c+1)>m,\qquad (b+1)+(c+1)>m.
\]
Hence
\[
        \cJ(\|I^{(m)}\|)=I^{(m-1)}\qquad (m\geq 2).
\]
Moreover,
\[
        \alpha(I^{(m)})=\left\lceil\frac{3m}{2}\right\rceil,
        \qquad
        \widehat{\alpha}(I)=\frac32 .
\]
The expected formula \eqref{eq:clean-formula} becomes
\[
        \pi_*\OO_{X'}\left(-\left(\left\lfloor\frac{3m}{2}\right\rfloor-2\right)F_0-(m-1)(F_1+F_2+F_3)\right).
\]
The \(F_i\)-conditions give exactly \(I^{(m-1)}\). The additional degree condition from \(F_0\) is automatic because
\[
        \alpha(I^{(m-1)})=\left\lceil\frac{3(m-1)}{2}\right\rceil\geq\left\lfloor\frac{3m}{2}\right\rfloor-2 .
\]
Thus the coordinate triangle has clean asymptotic multiplier ideals.
\end{nonitalicexample}

\begin{remark}
For reduced point sets in \(\PP^2\), the monomial situation is essentially exhausted by the preceding example and its one- or two-point subconfigurations. Indeed, if \(I=I(\mathcal P)\subset\CC[x,y,z]\) is a radical homogeneous monomial ideal defining a finite reduced subscheme of \(\PP^2\), then every minimal prime of \(I\) is one of
\[
        (x,y),\qquad (x,z),\qquad (y,z).
\]
Hence \(\mathcal P\) is a subset of the three coordinate points. In particular, an assumption that all symbolic powers \(I^{(m)}\), including \(I=I^{(1)}\), are monomial does not produce a genuinely broader class of reduced point configurations.
\end{remark}

\section{A conditional numerical obstruction for non-containment}

We now prove the conditional form of the multiplier-ideal obstruction.

\begin{proposition}\label{prop:clean-inclusion}
Assume that \(I=I(P_1,\ldots,P_s)\subset \CC[x_0,\ldots,x_n]\) has clean asymptotic multiplier ideals. Then, for every \(m\geq n\),
\[
        \cJ(\|I^{(m)}\|)
        \subseteq
        I^{(m+1-n)}.
\]
Moreover, equality holds if and only if
\begin{equation}\label{eq:equality-criterion}
        \lfloor m\widehat{\alpha}(I)\rfloor-n
        \leq
        \alpha(I^{(m+1-n)}).
\end{equation}
\end{proposition}

\begin{proof}
For integers \(a,b\geq 0\), the blow-up construction gives
\begin{equation}\label{eq:pushforward}
        \pi_*\OO_{X'}\Bigl(-aF_0-b\sum_{i=1}^sF_i\Bigr)
        =
        \bigoplus_{d\geq a}I^{(b)}_d .
\end{equation}
Indeed, \(F_0\) records the order at the cone vertex, which for a homogeneous form is its degree, while \(F_i\) records the order along the cone-line over \(P_i\).

Since \(I\) has clean asymptotic multiplier ideals, Definition~\ref{def:clean} and \eqref{eq:pushforward} give
\[
        \cJ(\|I^{(m)}\|)
        =
        \bigoplus_{d\geq \lfloor m\widehat{\alpha}(I)\rfloor-n}I^{(m+1-n)}_d .
\]
The inclusion in \(I^{(m+1-n)}\) is immediate. It is an equality exactly when the degree truncation is vacuous, namely exactly when \eqref{eq:equality-criterion} holds.
\end{proof}

Taking \(m=n\), we recover
\[
        \cJ(\|I^{(n)}\|)\subseteq I.
\]
The general Ein--Lazarsfeld--Smith argument gives this inclusion without the clean hypothesis. Here cleanliness is used to control precisely when the inclusion is strict, which is the point needed for the borderline power \(I^{(2r-1)}\).

We now specialize to \(n=2\).

\begin{theorem}\label{thm:conditional-main}
Let \(I=I(P_1,\ldots,P_s)\subset \CC[x,y,z]\) be the radical ideal of a finite set of points in \(\PP^2\). Assume that \(I\) has clean asymptotic multiplier ideals. If
\[
        I^{(2r-1)}\not\subseteq I^r,
\]
then
\begin{equation}\label{eq:main-bound}
        2r\widehat{\alpha}(I)\geq \alpha(I^{(2r-1)})+3.
\end{equation}
In particular,
\[
        \widehat{\alpha}(I)\geq 3.
\]
\end{theorem}

\begin{proof}
By Proposition~\ref{prop:clean-inclusion}, applied with \(m=2r\),
\[
        \cJ(\|I^{(2r)}\|)\subseteq I^{(2r-1)}.
\]
If equality held, then subadditivity of asymptotic multiplier ideals would give
\[
\begin{aligned}
        I^{(2r-1)}
        &=
        \cJ(\|I^{(2r)}\|) \\
        &\subseteq
        \cJ(\|I^{(2)}\|)^r \\
        &\subseteq
        I^r,
\end{aligned}
\]
contrary to the hypothesis. Thus the inclusion is strict. By the equality criterion \eqref{eq:equality-criterion},
\[
        \lfloor 2r\widehat{\alpha}(I)\rfloor-2>\alpha(I^{(2r-1)}).
\]
Since the right-hand side is an integer,
\[
        \lfloor 2r\widehat{\alpha}(I)\rfloor\geq \alpha(I^{(2r-1)})+3,
\]
and \eqref{eq:main-bound} follows.

Finally, by definition of \(\widehat{\alpha}(I)\) it holds
\[
        \widehat{\alpha}(I)\leq \frac{\alpha(I^{(2r-1)})}{2r-1}.
\]
Combining this with \eqref{eq:main-bound} yields
\[
        2r\widehat{\alpha}(I)\geq (2r-1)\widehat{\alpha}(I)+3,
\]
hence \(\widehat{\alpha}(I)\geq 3\).
\end{proof}

\begin{remark}
For \(r=2\), Theorem~\ref{thm:conditional-main} says that a clean ideal with \(I^{(3)}\not\subseteq I^2\) satisfies
\[
        4\widehat{\alpha}(I)\geq \alpha(I^{(3)})+3,
        \qquad\text{hence}\qquad
        \widehat{\alpha}(I)\geq\frac{\alpha_3+3}{4}.
\]
The general estimate
\[
        \frac{\alpha_m}{m+1}\leq \widehat{\alpha}(I)\leq\frac{\alpha_m}{m}
\]
(see \cite[Page~6]{Wal}) gives only \(\widehat{\alpha}(I)\geq\alpha_3/4\) from \(m=3\). Thus the term \(+3\) gives a genuine fixed-degree refinement. Combining the refined lower bound with \(\widehat{\alpha}(I)\leq\alpha_3/3\) recovers \(\widehat{\alpha}(I)\geq 3\); it does not strengthen the final assertion of Theorem~\ref{thm:conditional-main}.
\end{remark}

\section{The residual base locus problem}

Proposition~\ref{prop:dual-hesse-clean} shows why the multiplier-ideal strategy is not merely a computation of \(\widehat{\alpha}(I)\): the residual asymptotic base locus can be non-trivial even when it is ultimately invisible to the multiplier ideal. The natural problem is therefore valuation-theoretic.

\begin{problem}\label{prob:residual}
Let \(I=I(P_1,\ldots,P_s)\subset \CC[x,y,z]\) be the radical ideal of a finite set of points in \(\PP^2\). Describe the residual asymptotic base b-divisor \(\bR\) of Definition~\ref{def:residual-b-divisor}. Equivalently, by \eqref{eq:residual-order}, describe the asymptotic orders
\[
        v(\|I^{(m)}\|)
\]
for divisorial valuations \(v\) other than the expected valuations \(\ord_{F_0},\ord_{F_1},\ldots,\ord_{F_s}\). Determine, in particular, when all these extra valuations are redundant in the multiplier-ideal inequalities \eqref{eq:valuative-JM}, so that the clean formula \eqref{eq:clean-formula} holds.
\end{problem}

Before formulating a conjecture, we verify cleanliness for another family of line arrangements.

\begin{proposition}\label{prop:star-clean}
Let \(\mathcal A=L_1+\cdots+L_q\subset\PP^2\), \(q\geq 4\), be an arrangement of \(q\) lines with no three concurrent, and let
\[
        Z=\operatorname{Sing}(\mathcal A),\qquad I=I(Z).
\]
Then \(I\) has clean asymptotic multiplier ideals. Moreover,
\[
        \widehat{\alpha}(I)=\frac q2,
\]
and for every \(m\geq 2\),
\begin{equation}\label{eq:star-J}
        \cJ(\|I^{(m)}\|)=\bigoplus_{d\geq\lfloor mq/2\rfloor-2}I^{(m-1)}_d .
\end{equation}

The residual asymptotic base b-divisor \(\bR\) of Definition~\ref{def:residual-b-divisor} is non-zero. Its coefficient equals \(\frac12\) along the exceptional divisor obtained by blowing up the strict transform, in \(F_0\cong\operatorname{Bl}_Z\PP^2\), of each arrangement line. Here \(\frac12\) is the coefficient of that strict transform in the negative part of the Zariski decomposition of the class \(\frac q2H-\sum_{P\in Z}E_P\) on \(\operatorname{Bl}_Z\PP^2\).
\end{proposition}

\begin{proof}
Let \(s=\binom q2=|Z|\), and write
\[
        \rho:X=\operatorname{Bl}_Z\PP^2\longrightarrow\PP^2 .
\]
Denote by \(H\) the pull-back of a line, by \(E_P\) the exceptional divisor over a given \(P\in Z\), and by \(C_i\) the strict transform of \(L_i\). Then
\[
        C_i\sim H-\sum_{P\in L_i\cap Z}E_P,
        \qquad
        C_i^2=2-q,
\]
the \(C_i\)'s are pairwise disjoint, and
\[
        \sum_{i=1}^q C_i=qH-2\sum_{P\in Z}E_P .
\]
Set
\[
        B=(q-1)H-\sum_{P\in Z}E_P .
\]
If \(Q=L_1\cdots L_q\), then the degree \(q-1\) forms \(Q/L_i\) generate a linear system whose base locus in \(\PP^2\) is exactly \(Z\). At a node \(P=L_i\cap L_j\), the first-order parts of \(Q/L_i\) and \(Q/L_j\) are non-zero scalar multiples of local equations of \(L_j\) and \(L_i\), respectively. Hence they are transverse, and the strict transform of the system is base-point-free. Thus \(B\) is base-point-free and nef.

For any \(\lambda\in\mathbb R\) set
\[
        D_\lambda=\lambda H-\sum_{P\in Z}E_P .
\]
If \(q/2\leq\lambda\leq q-1\), then
\begin{equation}\label{eq:star-zariski}
        D_\lambda=\frac{2\lambda-q}{q-2}\,B+\frac{q-1-\lambda}{q-2}\sum_{i=1}^q C_i .
\end{equation}
Since \(B\cdot C_i=0\), the \(C_i\)'s are pairwise disjoint, and \(C_i^2=2-q<0\), \eqref{eq:star-zariski} is the Zariski decomposition. For \(\lambda\geq q-1\),
\[
        D_\lambda=B+(\lambda-q+1)H
\]
is base-point-free. In particular, if \(dH-k\sum E_P\) is effective, nefness of \(B\) gives
\[
        0\leq\Bigl(dH-k\sum E_P\Bigr)\cdot B=(q-1)d-\binom q2 k,
\]
so \(d/k\geq q/2\). Equality is attained by powers of \(Q\), because \(Q\) has degree \(q\) and multiplicity \(2\) at every point of \(Z\). Hence \(\widehat{\alpha}(I)=q/2\).

As in the proof of Proposition~\ref{prop:dual-hesse-clean}, after passing to the affine cone and blowing up the strict transforms \(C_i\subset F_0\), let \(G_i\) denote the new exceptional divisors. Let \(k=Nm\), let \(f\in I^{(k)}_d\) be homogeneous, and put \(\lambda=d/k\geq q/2\). If \(\lambda\leq q-1\), the coefficient of \(C_i\) in the negative part of \(kD_\lambda\) is
\[
        k\,\frac{q-1-\lambda}{q-2},
\]
while for \(\lambda\geq q-1\) there is no negative part. Thus
\[
        \frac1N v_{G_i}(f)\geq m\Bigl(\lambda+\max\Bigl\{0,\frac{q-1-\lambda}{q-2}\Bigr\}\Bigr).
\]
The expression in parentheses is non-decreasing in \(\lambda\), and at \(\lambda=q/2\) it equals \((q+1)/2\). Powers of \(Q\) show sharpness (choose \(Nm\) even). Therefore
\begin{equation}\label{eq:star-vG}
        v_{G_i}(\|I^{(m)}\|)=\frac{m(q+1)}{2}.
\end{equation}

Since \(v_{G_i}(F_0)=1\) and \(v_{G_i}(F_P)=0\) for all \(P\in Z\), \eqref{eq:residual-order} gives
\[
        r(\ord_{G_i})=\frac{q+1}{2}-\frac q2=\frac12>0,
\]
which is the coefficient of \(C_i\) in the negative part of \(D_{q/2}\). This proves \(\bR\neq0\). In particular, \(R_m\neq 0\) for every even \(m\): then \(\alpha_m=mq/2\) (attained by \(Q^{m/2}\)), so \(v_{G_i}(\frakb_m)\geq m(q+1)/2-mq/2=m/2>0\). For a single exponent \(m\), however, \(R_m\) may vanish along \(G_i\). For example, \(v_{G_i}(I)=q-1=\alpha_1\), the value being attained by \(Q/L_i\).

Again \(k_{G_i}=3\). Hence, besides the expected conditions
\[
        d\geq\Bigl\lfloor\frac{mq}{2}\Bigr\rfloor-2,
        \qquad
        f\in I^{(m-1)},
\]
the only possible extra condition is
\begin{equation}\label{eq:star-extra}
        v_{G_i}(f)\geq\Bigl\lfloor\frac{m(q+1)}{2}\Bigr\rfloor-3 .
\end{equation}

We show that \eqref{eq:star-extra} is automatic. Put
\[
        e_0=\Bigl\lfloor\frac{mq}{2}\Bigr\rfloor-2
\]
and let \(f\in I^{(m-1)}_d\) with \(d\geq e_0\). Fix \(L_i\) and write \(a=\ord_{L_i}(f)\).

If \(a\geq m-1\), the desired inequality \eqref{eq:star-extra} follows directly from \(v_{G_i}(f)=d+a\).

If otherwise \(a<m-1\), then B\'ezout (applied to \(L_i\) and the curve defined by \(f/L_i^a\)) gives
\[
        d-a\geq (q-1)(m-1-a)
\]
because \(L_i\) contains \(q-1\) points of \(Z\). Equivalently,
\begin{equation}\label{eq:star-bezout}
        a\geq\left\lceil\frac{(q-1)(m-1)-d}{q-2}\right\rceil,
\end{equation}
and thus
\begin{equation}\label{eq:star-lower}
        v_{G_i}(f)=d+a\geq d+\left\lceil\frac{(q-1)(m-1)-d}{q-2}\right\rceil .
\end{equation}
Since the right-hand side is non-decreasing in the integer variable \(d\), it attains its minimum for \(d=e_0\). If \(m=2t\), then
\[
        e_0=tq-2,\qquad (q-1)(m-1)-e_0=(t-1)(q-2)+1,
\]
so the right-hand side in \eqref{eq:star-bezout} equals \(t\), and \eqref{eq:star-lower} becomes
\[
        v_{G_i}(f)\geq t(q+1)-2=\frac{m(q+1)}{2}-2,
\]
which implies \eqref{eq:star-extra}. If \(m=2t+1\), then
\[
        e_0=tq+\Bigl\lfloor\frac q2\Bigr\rfloor-2,
\]
and
\[
        (q-1)(m-1)-e_0=t(q-2)-\Bigl\lfloor\frac q2\Bigr\rfloor+2 .
\]
For \(q\geq 4\), the right-hand side in \eqref{eq:star-bezout} is therefore also \(t\). Hence
\[
        v_{G_i}(f)\geq t(q+1)+\Bigl\lfloor\frac q2\Bigr\rfloor-2 .
\]
If \(q\) is even, this is one more than \(\lfloor m(q+1)/2\rfloor-3\); if \(q\) is odd, it is exactly \(\lfloor m(q+1)/2\rfloor-3\). Thus \eqref{eq:star-extra} holds in every case.

After removing the negative part in \eqref{eq:star-zariski}, the positive part is base-point-free, so there are no further residual valuations. The \(G_i\)-conditions are therefore redundant, and \eqref{eq:star-J} follows.
\end{proof}

Propositions~\ref{prop:dual-hesse-clean} and \ref{prop:star-clean} suggest the following conjecture.

\begin{conjecture}\label{conj:arrangements}
Let \(\mathcal A\) be a line arrangement in \(\PP^2\), and let
\[
        I=I(\operatorname{Sing}(\mathcal A)).
\]
Then \(I\) has clean asymptotic multiplier ideals.
\end{conjecture}

One could also conjecture a stronger statement, namely that \(\bR\) is always supported over the strict transforms of the arrangement lines in \(F_0\), i.e.\ that \(r(v)>0\) only for valuations whose centre on \(X'\) lies on one of these curves. This statement appears too rigid. For example, in the Klein and Wiman blow-ups, additional invariant curves of negative self-intersection, distinct from the arrangement divisor, are known (see \cite{Wal}). These curves may or may not contribute to \(\bR\), and even if they do, they might still impose redundant conditions for the asymptotic multiplier ideal. Thus the clean formulation separates the property actually needed for Theorem~\ref{thm:conditional-main} from a much stronger assertion about the complete base locus.

\begin{question}
For a line arrangement \(\mathcal A\), to what extent are the residual asymptotic orders of \(I(\operatorname{Sing}(\mathcal A))\), i.e.\ the b-divisor \(\bR\), and in particular the ideals \(\cJ(\|I^{(m)}\|)\), determined by the incidence combinatorics of \(\mathcal A\)?
\end{question}

The next natural test cases are the Fermat--Ceva, Klein, and Wiman arrangements. The latter two are especially informative because their blow-ups already carry non-arrangement negative curves \cite{Wal}.

There is also a substantially stronger question within the present point-ideal framework.

\begin{question}
Does every radical ideal of a finite reduced set of points in \(\PP^2\) have clean asymptotic multiplier ideals?
\end{question}

\begin{remark}
The residual-base-locus problem also clarifies why the multiplier-ideal method cannot be expected to yield, without further hypotheses, a purely numerical criterion for the containment \(I^{(2r-1)}\subseteq I^r\) in terms of the invariants \(\alpha_m\) and \(\widehat{\alpha}(I)\) alone. An example of such a criterion is the contrapositive of Theorem~\ref{thm:conditional-main}: if \(2r\widehat{\alpha}(I)<\alpha(I^{(2r-1)})+3\), then \(I^{(2r-1)}\subseteq I^r\). Multiplier ideals are sensitive to the full asymptotic base locus, encoded in \(\bR\), not only to the smallest degree in which a symbolic power is non-zero. Cleanliness is precisely the hypothesis that \(\bR\) is invisible to them. This limitation concerns the multiplier-ideal approach only: sufficient conditions for containment involving other invariants, such as the Castelnuovo--Mumford regularity, are available; see \cite{BH}.
\end{remark}

\section*{Acknowledgments}
The road to this paper has been a long one: the work began in 2016 and, at times, followed a winding path. We are confident that the effort has proven worthwhile. We are sincerely grateful to Leibniz Universit\"at Hannover for hosting us. Our time there left us with many fond memories and provided the collaborative environment from which this work emerged.
\section*{Funding}
The research was partly conducted in the framework of the DFG-funded research training group RTG 2965: From Geometry to Numbers (Project number 512730679).

\section*{AI Declaration}
Back when we began this project, generative AI wasn’t even on the map.

\bigskip

\noindent
V\'ictor Gonz\'alez-Alonso,\\
Institut f\"ur Algebraische Geometrie,\\
Leibniz Universit\"at Hannover,\\
Welfengarten 1, D-30167 Hannover, Germany.\\
\textit{E-mail address:} \texttt{gonzalez@math.uni-hannover.de}

\bigskip

\noindent
Piotr Pokora,\\
Department of Mathematics,\\
University of the National Education Commission Krakow\\
Podchor\c a\.zych 2, PL-30-084 Krak\'ow, Poland.\\
\textit{E-mail address:} \texttt{piotr.pokora@uken.krakow.pl}


\begin{thebibliography}{99}\footnotesize

\bibitem{BNC}
Th. Bauer, S. Di Rocco, B. Harbourne, J. Huizenga, A. Lundman, P. Pokora, and T. Szemberg:
Bounded negativity and arrangements of lines.
\textit{Int. Math. Res. Notices} \textbf{2015(19)}: 9456--9471 (2015).

\bibitem{Wal}
Th. Bauer, S. Di Rocco, B. Harbourne, J. Huizenga, A. Seceleanu, and T. Szemberg:
Negative curves on symmetric blowups of the projective plane, resurgences and Waldschmidt constants.
\textit{Int. Math. Res. Notices} \textbf{2019(24)}: 7459--7514 (2019).

\bibitem{BH}
C. Bocci and B. Harbourne:
Comparing powers and symbolic powers of ideals.
\textit{J. Algebraic Geometry} \textbf{19}: 399--417 (2010).

\bibitem{DST}
M. Dumnicki, T. Szemberg, and H. Tutaj-Gasi\'nska:
Counterexamples to the \(I^{(3)} \subset I^2\) containment.
\textit{J. Algebra} \textbf{393}: 24--29 (2013).

\bibitem{ELS}
L. Ein, R. Lazarsfeld, and K. Smith:
Uniform bounds and symbolic powers on smooth varieties.
\textit{Invent. Math.} \textbf{144}: 241--252 (2001).

\bibitem{HaHu}
B. Harbourne and C. Huneke:
Are symbolic powers highly evolved?
\textit{J. Ramanujan Math. Soc.} \textbf{28A}: 247--266 (2013).

\bibitem{HH}
M. Hochster and C. Huneke:
Comparison of symbolic and ordinary powers of ideals.
\textit{Invent. Math.} \textbf{147}: 349--369 (2002).

\bibitem{HowaldMonomialMultiplierIdeals}
J.~A. Howald:
Multiplier ideals of monomial ideals.
\textit{Trans. Amer. Math. Soc.} \textbf{353}(7): 2665--2671 (2001).


\bibitem{JM}
M. Jonsson and M. Musta\c{t}\u{a}:
Valuations and asymptotic invariants for sequences of ideals.
\textit{Ann. Inst. Fourier (Grenoble)} \textbf{62}(6): 2145--2209 (2012).

\bibitem{Klein}
F. Klein:
\"Uber die Transformation siebenter Ordnung der elliptischen Functionen.
\textit{Math. Annalen} \textbf{14}: 428--471 (1879).


\bibitem{KM}
J. Koll\'ar and S. Mori:
\emph{Birational Geometry of Algebraic Varieties}.
Cambridge Tracts in Mathematics, vol.~134, Cambridge University Press, Cambridge, 1998.

\bibitem{LazarsfeldPositivityII}
R.~Lazarsfeld:
\emph{Positivity in Algebraic Geometry II. Positivity for Vector Bundles, and Multiplier Ideals},
Ergebnisse der Mathematik und ihrer Grenzgebiete, 3. Folge, vol.~49,
Springer-Verlag, Berlin, 2004.

\bibitem{MSS}
G. Malara, T. Szemberg, and J. Szpond:
On a conjecture of Demailly and new bounds on Waldschmidt constants in \(\PP^N\).
\textit{J. Number Theory} \textbf{189}: 211--219 (2018).

\bibitem{Sec}
A. Seceleanu:
A homological criterion for the containment between symbolic and ordinary powers of some ideals of points in \(\PP^2\).
\textit{J. Pure Appl. Algebra} \textbf{219}: 4857--4871 (2015).

\bibitem{SzSz}
T. Szemberg and J. Szpond:
On the containment problem.
\textit{Rend. Circ. Mat. Palermo (2)} \textbf{66(2)}: 233--245 (2017).

\bibitem{Teitler}
Z. Teitler:
Bounding symbolic powers via asymptotic multiplier ideals.
\textit{Ann. Univ. Paedagog. Crac. Stud. Math.} \textbf{8}: 67--77 (2009).

\bibitem{Wiman}
A. Wiman:
Zur Theorie der endlichen Gruppen von birationalen Transformationen in der Ebene.
\textit{Math. Annalen} \textbf{48}: 195--240 (1896).

\end{thebibliography}
\end{document}